# STABILITY FOR HOLOMORPHIC SPHERES AND MORSE THEORY

RALPH L. COHEN, JOHN D.S. JONES, AND GRAEME B. SEGAL

ABSTRACT. In this paper we study the question of when does a closed, simply connected, integral symplectic manifold $(X, \omega)$ have the stability property for its spaces of based holomorphic spheres? This property states that in a stable limit under certain gluing operations, the space of based holomorphic maps from a sphere to $X$, becomes homotopy equivalent to the space of all continuous maps,
$$\varinjlim Hol_{x_0}(\mathbb{P}^1, X) \simeq \Omega^2 X.$$
This limit will be viewed as a kind of *stabilization* of $Hol_{x_0}(\mathbb{P}^1, X)$. We conjecture that this stability property holds if and only if an evaluation map $E : \varinjlim Hol_{x_0}(\mathbb{P}^1, X) \to X$ is a quasifibration. In this paper we will prove that in the presence of this quasifibration condition, then the stability property holds if and only if the Morse theoretic flow category (defined in [4]) of the symplectic action functional on the $\mathbb{Z}$ - cover of the loop space, $\tilde{L}X$, defined by the symplectic form, has a classifying space that realizes the homotopy type of $\tilde{L}X$. We conjecture that in the presence of this quasifibration condition, this Morse theoretic condition always holds. We will prove this in the case of $X$ a homogeneous space, thereby giving an alternate proof of the stability theorem for holomorphic spheres for a projective homogeneous variety originally due to Gravesen [7].

## INTRODUCTION

Let $(X, \omega, J)$ be a closed, connected, integral symplectic manifold of dimension $2n$ with a compatible almost complex structure. Here $\omega$ is the symplectic 2 - form and $J$ is the almost complex structure. By "integral" symplectic manifold we mean that the symplectic form $\omega$ defines an *integer* cohomology class, $[\omega] \in H^2(X;\mathbb{Z})$. Recall that a map $f : \mathbb{CP}^1 \to X$ is $J$ - holomorphic if $df \circ j = J \circ df$, where $j$ is the almost complex structure on the tangent bundle of $\mathbb{CP}^1$. Let $x_0 \in X$ be a fixed based point. In this paper we consider the space of based $J$ - holomorphic spheres,

$$Hol_{x_0}(\mathbb{P}^1, X) = \{f : \mathbb{CP}^1 \to X, \text{such that } f \text{ is } J \text{ - holomorphic and } f(0) = x_0.\}$$

where $0 \in \mathbb{C} \subset \mathbb{C} \cup \infty$ is taken to be the basepoint of $\mathbb{CP}^1$. The holomorphic mapping space $Hol_{x_0}(\mathbb{P}^1, X)$ is topologized as a subspace of the two fold loop space, $\Omega^2 X$. The relative homotopy type of the pair $Hol_{x_0}(\mathbb{P}^1, X) \subset \Omega^2 X$ has been studied for a variety of complex manifolds $X$. For example for $X = \mathbb{CP}^n$, it was proved in [20] that

$$\varinjlim_k Hol^k_{x_0}(\mathbb{CP}^1, \mathbb{CP}^n) \simeq \Omega^2 \mathbb{CP}^n$$

Date: April 29, 1999.
The first author was partially supported by a grant from the NSF.
The second author was partially supported by the American Institute of Mathematics.

where the subscript $k$ denotes the degree (or the homology class) of the mapping, and the limit actually refers to a homotopy colimit of spaces under the gluing of a fixed map $\iota \in Hol^1_{x_0}(\mathbb{CP}^1, \mathbb{CP}^n)$. Similar stability theorems have been proven for $X$ a Grassmannian or more general flag manifold [8], [11], [7] $X$ a toric variety [9], and for $X$ a loop group, $X = \Omega G$, where $G$ is a compact Lie group [7],[21]. The purpose of this paper is to begin a general investigation of what basic properties of the symplectic manifold $(X, \omega)$ assure that there is an appropriate limiting process so that

$$\varinjlim Hol_{x_0}(\mathbb{P}^1, X) \simeq \Omega^2 X.$$

We will refer to this property as the *stability property* for the space of $J$ - holomorphic spheres in $(X, \omega)$.

To make this question more precise, in section 1 we describe a space $Hol_{x_0}(\mathbb{P}^1, X)^+$, built out of limits of "chains" of holomorphic maps, which is an appropriate "stabilization" of the holomorphic mapping space $Hol_{x_0}(\mathbb{P}^1, X)$. To do this we will assume the following positivity condition.

**Definition 1.** We will say that the symplectic manifold $(X, \omega)$ is *positive*, if there is a real number $\lambda > 0$ such that

$$\langle c_1(X), \gamma \rangle \geq \lambda \langle \omega, \gamma \rangle$$

for all $\gamma \in \pi_2(X)$.

The stabilization $Hol_{x_0}(\mathbb{P}^1, X)^+$ of $Hol_{x_0}(\mathbb{P}^1, X)$, which we will construct in the next section will have the the following properties.

1. The inclusion of holomorphic maps into continuous maps $Hol_{x_0}(\mathbb{P}^1, X) \hookrightarrow \Omega^2 X$ naturally extends to a map

$$j : Hol_{x_0}(\mathbb{P}^1, X)^+ \to \Omega^2 X.$$

2. Suppose that there are gluing operations

$$Hol^k_{x_0}(\mathbb{CP}^1, X) \times Hol^r_{x_0}(\mathbb{CP}^1, X) \to Hol^{k+r}_{x_0}(\mathbb{CP}^1, X)$$

that lift (up to homotopy) the loop sum operation

$$\Omega^2 X \times \Omega^2 X \to \Omega^2 X.$$

Here $Hol^q$ denotes the subspace of $Hol_{x_0}(\mathbb{P}^1, X)$ consisting of classes represented by maps $\theta : \mathbb{P}^1 \to X$ with $\langle [\omega], [\theta] \rangle = k \in \mathbb{Z}$, where $[\omega] \in H^2(X, \mathbb{Z})$ is the class represented by the symplectic form. Let $\iota \in Hol^q$ be any fixed class with $q \neq 0$. Then consider the map

$$*\iota : Hol^k_{x_0}(\mathbb{CP}^1, X) \to Hol^{k+q}_{x_0}(\mathbb{CP}^1, X)$$

given by gluing with $\iota$. Then $Hol_{x_0}(\mathbb{P}^1, X)^+$ is homotopy equivalent to the homotopy colimit of this operation:

$$Hol_{x_0}(\mathbb{P}^1, X)^+ \simeq \varinjlim{}_k Hol^k_{x_0}(\mathbb{CP}^1, X).$$

Furthermore, if the gluing operation on $Hol_{x_0}(\mathbb{P}^1, X)$ gives it the structure of a $\mathcal{C}_2$- operad space in the sense of May [14], and $\pi_0(Hol_{x_0}(\mathbb{P}^1, X))$ is finitely generated as a monoid, then $Hol_{x_0}(\mathbb{P}^1, X)^+$ is the group completion

$$Hol_{x_0}(\mathbb{P}^1, X)^+ \simeq \Omega B(Hol_{x_0}(\mathbb{P}^1, X)).$$

The example to think of in this case is $X$ a Grassmannian or more general flag manifold.

**Note:** Since a two fold loop space is a $\mathcal{C}_2$ -operad space, one can take its group completion, $\Omega^2 X^+ = \Omega B(\Omega^2 X)$. But clearly $\Omega^2 X^+ = \Omega^2 X$, and so in this case, when $Hol_{x_0}(\mathbb{P}^1, X)$ has the structure of a $\mathcal{C}_2$ - space, then the map $j : Hol_{x_0}(\mathbb{P}^1, X)^+ \to \Omega^2 X$ mentioned in property 1, is simply given by applying the group completion functor to the inclusion $Hol_{x_0}(\mathbb{P}^1, X) \hookrightarrow \Omega^2 X$.

We can now phrase the question posed above more precisely: What geometric properties of a symplectic manifold $(X, \omega)$, assure that $Hol_{x_0}(\mathbb{P}^1, X)^+ \simeq \Omega^2 X$? We conjecture an answer to this question as follows.

Consider the evaluation map

$$E : Hol_{x_0}(\mathbb{P}^1, X) \to X$$
$$f \to f(\infty)$$

As we will see, this evaluation map naturally extends to the group completion

$$E : Hol_{x_0}(\mathbb{P}^1, X)^+ \to X.$$

This extension can be viewed as the composition $Hol_{x_0}(\mathbb{P}^1, X)^+ \xrightarrow{j} \Omega^2 X \xrightarrow{E} X$. Recall that the basepoint condition is that $f(0) = x_0$. The map $E$ evaluates a map at the other pole, $\infty \in \mathbb{C} \cup \infty = \mathbb{CP}^1$. Recall that $X$ is said to be *rationally connected* if each of the fibers, $Hol_{x_0, x}(\mathbb{CP}^1, X) = E^{-1}(x) \subset Hol_{x_0}(\mathbb{P}^1, X)$ is nonempty. A stronger condition is that the map $E : Hol_{x_0}(\mathbb{P}^1, X)^+ \to X$ is a quasifibration. This in particular implies that all of its fibers $Hol_{x_0, x}(\mathbb{CP}^1, X)^+ = E^{-1}(x) \subset Hol_{x_0}(\mathbb{P}^1, X)^+$ are homotopy equivalent.

**Definition 2.** We say that a symplectic manifold $(X, \omega)$ has the "quasifibration property" if the evaluation map

$$E : Hol_{x_0}(\mathbb{P}^1, X)^+ \to X$$

is a quasifibration.

**Conjecture.** *A symplectic manifold has the quasifibration property if and only if*

$$Hol_{x_0}(\mathbb{P}^1, X)^+ \simeq \Omega^2 X.$$

In this paper we address this conjecture using Morse theory techniques. In particular we use the approach to Morse theory introduced by the authors in [4], [5], [3], which studies the *flow category*, $\mathcal{C}_f$, of a function $f : M \to \mathbb{R}$, where $M$ is a smooth manifold with no boundary. Specifically, $\mathcal{C}_f$ is a topological category (where both the objects and morphisms are topologized), with the objects of $\mathcal{C}_f$ being the critical points of $f$, and the morphisms between a critical point $a$ and $b$, $Mor_{\mathcal{C}_f}(a,b)$ are given by the space of *piecewise flow lines*, $\bar{\mathcal{M}}(a,b)$. If $M$ is a closed, finite dimensional manifold, then in a generic setting, (i.e $f : M \to \mathbb{R}$ a Morse function satisfying the Morse - Smale transversality condition) the space of flow lines $\mathcal{M}(a,b)$ is a smooth, framed manifold of dimension $\lambda(a) - \lambda(b) - 1$, where $\lambda$ refers to the index of the critical point. $\mathcal{M}(a,b)$ is not compact, however, but is, when we adjoin all piecewise flow lines

$$\bar{\mathcal{M}}(a,b) = \bigcup \mathcal{M}(a,a_1) \times \mathcal{M}(a_1,a_2) \times \cdots \times \mathcal{M}(a_r,b)$$

where the union is taken over all sets of critical points $a_1, \cdots, a_r$, such that the above flow spaces are nonempty. The details of this flow category, including the topology of the morphism spaces, will be reviewed in section 1. Let $B\mathcal{C}_f$ be the classifying space (that is the geometric realization of the simplicial nerve) of the category $\mathcal{C}_f$. The main theorem of [4] was the following.

**Theorem 1.** *Let $f : M \to \mathbb{R}$ be a Morse function on a closed finite dimensional Riemannian manifold $M$. Then there is a natural map $\phi : B\mathcal{C}_f \to M$ satifying the following:*

1. *If the Smale transversality condition is satisfied, then $\phi$ is a homeomorphism.*
2. *Even without the transversality condition $\phi$ is a homotopy equivalence.*

In this paper we will address the above conjecture on the homotopy type of $Hol_{x_0}(\mathbb{P}^1, X)^+$ by studying the flow category of the symplectic area functional on the loop space of a symplectic manifold. More specifically, let $(X, \omega)$ be a simply connected closed, finite dimensional symplectic manifold, and let $LX$ denote the space of smooth loops,

$$LX = C^\infty(S^1, X).$$

Let $\tilde{L}X$ be the $\mathbb{Z}$ - cover of the loop space defined by the homomorphism

(0.1) $$\pi_1(LX) = \pi_2(X) \to H_2(X) \xrightarrow{\omega} \mathbb{Z}$$

$$\tilde{L}X = \{(\gamma, \theta) \in LX \times C^\infty(D^2, X) : \theta_{|S^1} = \gamma\}/\sim$$

where the equivalence relation is given by $(\gamma, \theta_1) \sim (\gamma, \theta_2)$, if the map $\theta : S^2 \to X$ defined to be $\theta_1$ on the upper hemisphere, and $\theta_2$ on the lower hemisphere, has the property that $\omega([\theta]) = 0 \in \mathbb{Z}$.

The symplectic area function is the map

$$\alpha : \tilde{L}X \to \mathbb{R}$$

$$(\gamma, \theta) \to \int_{D^2} \theta^*(\omega).$$

This is the functional that is the basis of symplectic Floer homology [6]. The critical points are pairs $(x, n) \in X \times \mathbb{Z}$, where a point $x \in X$ is viewed as a constant loop extended to a map $\theta : (D^2, S^1) \to (X, x)$ with $\omega([\theta]) = n$. The critical points are topologized as $X \times \mathbb{Z}$. Moreover a gradient flow line of $\alpha$ from a critical point $(x_1, n_1)$ to a critical point $(x_2, n_2)$ is a $J$ - holomorphic map $\gamma \in Hol_{x_1, x_2}(\mathbb{CP}^1, X)$ that represents an element in $\pi_2(X)$ with $\omega([\gamma]) = n_1 - n_2 \in \mathbb{Z}$. Let $\mathcal{C}_\alpha$ be the flow category of the symplectic area function. Notice that the above theorem does not immediately apply to the classifying space $B\mathcal{C}_\alpha$ because the loop space $LX$ is an infinite dimensional manifold. Our main result is that this is essentially the main issue in attacking the above conjecture. Namely, we will prove the following theorem.

**Theorem 2.** *Let $(X, \omega)$ be a simply connected, positive, integral symplectic manifold that has the quasifibration property. Then the map*

$$j : Hol_{x_0}(\mathbb{P}^1, X)^+ \to \Omega^2 X$$

*is a homotopy equivalence if and only if the map*

$$B\mathcal{C}_\alpha \to \tilde{L}X$$

*is a homotopy equivalence.*

The paper is organized as follows. In section one we recall in more detail the construction of the classifying space of a flow category and then define and study the stabilization $Hol_{x_0}(\mathbb{P}^1, X)^+$. In section 2 we prove theorem 2. In section 3 we use this theorem to give an alternate proof of a theorem of Gravesen [7], stating that for $X$ a homogeneous space, its space of $J$ - holomorphic curves has the stability property. We will continue our study of the flow category $\mathcal{C}_\alpha$ and therefore of the above conjecture in future work.

1. Classifying spaces and stabilized holomorphic spheres

In this section we describe the classifying space for the flow category of the symplectic action functional, and the stabilization $Hol_{x_0}(\mathbb{P}^1, X)^+$ discussed in the introduction.

Given a smooth function $f : M \to \mathbb{R}$ then recall from [4] that the flow category $\mathcal{C}_f$ has objects the critical points of $f$. As mentioned in the introduction, the morphism space is given by the space of piecewise flow lines, whose definition (including its topology) we recall now.

Let $\phi : \mathbb{R} \to M$ be a smooth curve in $M$. The *flow equation* is the differential equation

(1.1) $$\frac{d\phi}{dt} = \nabla f(\phi(t)).$$

The curves which satisfy this equation are the *flow lines* of $f$.

Let $\phi$ be a flow line starting at $a$ and ending at $b$. Since $f$ is strictly increasing along flow-lines it defines a diffeomorphism of the flow-line $\phi$ with the open interval $(f(a), f(b))$. The inverse of this diffeomorphism gives a parametrization of the flow-line as a smooth curve

$$\gamma : (f(a), f(b)) \to M$$

such that

$$f(\gamma(t)) = t.$$

Furthermore, $\gamma$ extends to a continuous function $\gamma : [f(a), f(b)] \to M$, which we will call an *extended flow-line*, by setting $\gamma(f(a)) = a$ and $\gamma(f(b)) = b$. We will refer to this parametrisation of a flow-line, or an extended flow-line, as the *parametrization by level sets*. A flow-line parametrized by level sets satisfies the differential equation

(1.2) $$\frac{d\gamma}{dt} = \frac{\nabla f(\phi(t))}{\|\nabla f(\phi(t))\|^2}.$$

Define $\bar{\mathcal{M}}(a, b) = \bar{\mathcal{M}}_f(a, b)$ to be the space of continuous curves $\gamma : [f(a), f(b)] \to M$ with the following properties.

1. The curve $\gamma(t)$ passes through at most a finite number of critical points of $f$.
2. On the complement of the set of $t$ where $\gamma(t)$ is a critical point, the curve $\gamma$ is smooth and satisfies equation 1.2
3. The curve $\gamma$ starts at $a$ and ends at $b$, that is $\gamma(f(a)) = a$ and $\gamma(f(b)) = b$.

This space $\bar{\mathcal{M}}(a, b)$ is topologized as a subspace of the space $C^0([f(a), f(b)], M)$ of continuous functions $[f(a), f(b)] \to M$ with the topology of uniform convergence. A *piecewise flow-line* is a curve which satisfies these three properties and $\bar{\mathcal{M}}(a, b)$ is the space of piecewise flow-lines from $a$ to $b$. If $\gamma$ is in $\bar{\mathcal{M}}(a, b)$ then when we remove the finite set of points where $\gamma(t)$ is a critical point of $f$ then each component of the resulting curve is a flow-line parametrized by level sets. Thus a curve in $\bar{\mathcal{M}}(a, b)$ is a piecewise flow-line from $a$ to $b$ in the informal sense described in the introduction.

Notice that there are natural pairing operations

$$\bar{\mathcal{M}}(a, b) \times \bar{\mathcal{M}}(b, c) \hookrightarrow \bar{\mathcal{M}}(a, c)$$

$$\gamma_1 \times \gamma_2 \to \gamma_1 * \gamma_2$$

where $\gamma_1 * \gamma_2$ is simply the concantenation of piecewise flows. This operation makes $\mathcal{C}_f$ into a category, where we note that the identity morphism $1 \in \bar{\mathcal{M}}(a, a)$ is the constant flow at $a$.

In the case when $f : M \to \mathbb{R}$ is a Morse function, a specific map

$$\phi : B\mathcal{C}_f \to M$$

was described, and was shown to be a homeomorphism when the Morse - Smale transversality conditions are satisfied, and a homotopy equivalence in general. The homotopy type of the map $\phi$ has a rather easy description, which we now give.

Given a space $X$ Let $\mathcal{C}(X)$ be the category whose objects are the points of $X$, and whose morphisms $Mor(x_1, x_2)$ are all continuous paths defined on some closed interval, that begin at $x_1$ and end at $x_2$.

$$Mor(x_1, x_2) = \{\beta : [a, b] \to M : \beta(a) = x_1, \ \beta(b) = x_2, \text{ for some interval } [a, b]\}.$$

It is a standard fact that the classifying space $B\mathcal{C}(X)$ is natually homotopy equivalent to $X$ (see [18]).

Notice that for a smooth function $f : M \to \mathbb{R}$, the flow category $\mathcal{C}_f$ is a subcategory of $\mathcal{C}(M)$. The inclusion $\iota : \mathcal{C}_f \hookrightarrow \mathcal{C}(M)$ induces on the classifying space level the map mentioned above.

**Definition 3.** Given a smooth function $f : M \to \mathbb{R}$, define $\phi : B\mathcal{C}_f \to M$ to be the composition

$$\phi : B\mathcal{C}_f \xrightarrow{B\iota} B\mathcal{C}(M) \simeq M.$$

We now consider the symplectic action functional as described in the introduction. So let $(X, \omega)$ be a closed, simply connected, integral symplectic manifold, and let $LX$ and $\tilde{L}X$ denote the loop space and its $\mathbb{Z}$ - cover as defined in the introduction. Also as in the introduction, let

$$\alpha : \tilde{L}X \to \mathbb{R}$$

denote the symplectic action functional, and $\mathcal{C}_\alpha$ its flow category. The objects in $\mathcal{C}_\alpha$ are the critical points of $\alpha$, which are given by pairs $(x, n) \in X \times \mathbb{Z}$. Thus the objects of $\mathcal{C}_\alpha$ have a nontrivial topology. Recall that $(x, n) \in Obj(\mathcal{C}_\alpha)$ corresponds to the constant loop at $x \in X$, extended to a map of a disk $\theta : (D^2, S^1) \to (X, x)$ so that $\omega([\theta]) = n \in \mathbb{Z}$. As mentioned in the introduction, a flow from $(x_1, n_1)$ to $(x_2, n_2)$ of $\alpha$ is given by a holomorphic map $\phi \in Hol^{n_1-n_2}_{x_1, x_2}(\mathbb{P}^1, X)$, where the subscript denotes the image under the map $\phi$ of the poles $0$ and $\infty$ in $\mathbb{P}^1$, and the superscript denotes the value $\omega([\phi]) \in \mathbb{Z}$. (See Floer's original paper [6] for details on the dynamics of the symplectic action map $\alpha$.) Thus the morphism space $Mor_{\mathcal{C}_\alpha}((x_1, n_1); (x_2, n_2))$ is given by the space of *piecewise flows* whose topology is as described above. We think of this space as the space of *piecewise holomorphic maps* which we denote $\bar{Hol}^{n_1-n_2}_{x_1, x_2}(\mathbb{P}^1, X)$. An element of this space can be viewed as a chain of holomorphic maps

$$\phi = \phi_1 \vee \phi_2 \vee \cdots \vee \phi_k$$

where each $\phi_i : \mathbb{P}^1 \to X$ is a holomorphic map satisfying the following:

1. $\phi_1(0) = x_1$,

2. $\phi_i(\infty) = \phi_{i+1}(0) \quad \text{for } i = 1, \cdots, k - 1$

3. $\phi_k(\infty) = x_2$

4. The homotopy class represented by the composition

$$\phi : S^2 \xrightarrow{fold} \bigvee_k S^2 \xrightarrow{\phi_1 \vee \cdots \vee \phi_k} X$$

has the property that $\omega([\phi]) = n_1 - n_2 \in \mathbb{Z}$.

The space $\bar{Hol}^{n_1-n_2}_{x_1,x_2}(\mathbb{P}^1, X)$ can be viewed as a partial compactification of the holomorphic mapping space $Hol^{n_1-n_2}_{x_1,x_2}(\mathbb{P}^1, X)$, and in particular it maps to a subspace of the space of stable curves as described by Kontsevich and Manin [13]. An important difference in the topology of $\bar{Hol}_{x_1,x_2}(\mathbb{P}^1, X)$ and that of the Kontsevich - Manin moduli space of stable curves is that they consider the orbits of holomorphic maps of spheres under the action of the (holomorphic) automorphisms of $\mathbb{P}^1$. We do not divide out by this group of parameterizations.

Let $\bar{Hol}^n_{x_0}(\mathbb{P}^1, X)$ denote the union of the spaces $\bar{Hol}^n_{x_0,y}(\mathbb{P}^1, X)$ as $y \in X$ varies. It is topologized in a natural way so that the evaluation map

(1.3) $$E : \bar{Hol}^n_{x_0}(\mathbb{P}^1, X) \to X$$
$$\phi = \phi_1 \vee \cdots \vee \phi_k \to \phi_k(\infty)$$

is continuous. Notice that we have a continuous inclusion

$$\bar{Hol}^k_{x_0}(\mathbb{P}^1, X) \hookrightarrow \Omega^2 X$$

whose image lies in the components of $\Omega^2 X$ representing homotopy classes which map to $k \in \mathbb{Z}$ under the homomorphism

$$\omega : \pi_2(X) \to \mathbb{Z}.$$

Notice that there is a monoid structure on $\bar{Hol}^*_{x_0,x_0}(\mathbb{P}^1, X)$

(1.4) $$\bar{Hol}^{n_1}_{x_0,x_0}(\mathbb{P}^1, X) \times \bar{Hol}^{n_2}_{x_0,x_0}(\mathbb{P}^1, X) \to \bar{Hol}^{n_1+n_2}_{x_0,x_0}(\mathbb{P}^1, X)$$

given by concantenations of piecewise holomorphic maps:

$$(\gamma = \gamma_1 \vee \cdots \vee \gamma_r) \times (\phi = \phi_1 \vee \cdots \vee \phi_k) \to \gamma_1 \vee \cdots \vee \gamma_r \vee \phi_1 \vee \cdots \vee \phi_k.$$

This also extends to give a natural action

(1.5) $$\bar{Hol}^{n_1}_{x_0,x_0}(\mathbb{P}^1, X) \times \bar{Hol}^{n_2}_{x_0}(\mathbb{P}^1, X) \to Hol^{n_1+n_2}_{x_0}(\mathbb{P}^1, X).$$

We will now use this action to define the stabilization $Hol_{x_0}(\mathbb{P}^1, X)^+$ of $Hol_{x_0}(\mathbb{P}^1, X)$. As mentioned in the introduction, $Hol_{x_0}(\mathbb{P}^1, X)^+$ will be a certain limit of holomorphic mapping spaces. We now make this idea precise.

Let $(X, \omega)$ satisfy the positivity property 1. Choose a fixed $\gamma \in \bar{Hol}^n_{x_0,x_0}(\mathbb{P}^1, X)$ with $n \neq 0$. Consider the map

$$\gamma * : \bar{Hol}_{x_0}(\mathbb{P}^1, X) \to \bar{Hol}_{x_0}(\mathbb{P}^1, X)$$

given by acting by $\gamma$ as in 1.5. We define $Hol_{x_0}(\mathbb{P}^1, X)^+$ to be the homotopy colimit under this gluing map:

**Definition 4.** The "stabilization" $Hol_{x_0}(\mathbb{P}^1, X)^+$ of the holomorphic mapping space is the homotopy colimit

$$Hol_{x_0}(\mathbb{P}^1, X)^+ = hoco\varinjlim \bar{Hol}_{x_0}(\mathbb{P}^1, X)$$

where the the homotopy colimit is taken with respect to the gluing map $\gamma * : \bar{Hol}_{x_0}(\mathbb{P}^1, X) \to \bar{Hol}_{x_0}(\mathbb{P}^1, X)$.

We observe a couple of properties of this construction. First, notice that the following diagram homotopy commutes

$$\begin{array}{ccc} \bar{Hol}_{x_0}(\mathbb{P}^1, X) & \xrightarrow{\gamma *} & \bar{Hol}_{x_0}(\mathbb{P}^1, X) \\ \cap \downarrow & & \downarrow \cap \\ \Omega^2 X & \xrightarrow{\gamma *} & \Omega^2 X \end{array}$$

where the bottom horizontal map represents the "loop sum" operation with the element $\gamma$. Notice furthermore that the map $\gamma : \Omega^2 X \to \Omega^2 X$ is a homotopy equivalence, with homotopy inverse given by taking the loop sum operation with an element $\gamma^{-1} \in \Omega^2 X$ representing $-[\gamma] \in \pi_2(X)$. Thus the inclusion map $\bar{Hol}_{x_0}(\mathbb{P}^1, X) \hookrightarrow \Omega^2 X$ extends to a map

(1.6) $$j : Hol_{x_0}(\mathbb{P}^1, X)^+ \hookrightarrow \Omega^2 X.$$

Notice furthermore that the evaluation map ( 1.3) $E : \bar{Hol}_{x_0}(\mathbb{P}^1, X) \to X$ commutes with the left action by $\gamma$, in the sense that the following diagram commutes:

$$\begin{array}{ccc} \bar{Hol}_{x_0}(\mathbb{P}^1, X) & \xrightarrow{\gamma *} & \bar{Hol}_{x_0}(\mathbb{P}^1, X) \\ E \downarrow & & \downarrow E \\ X & \xrightarrow{=} & X \end{array}$$

Thus the evaluation map extends to a map of the group completion

(1.7) $$E : Hol_{x_0}(\mathbb{P}^1, X)^+ \to X.$$

With these definitions (of the stabilization and the evaluation map), we may now recall from definition 2 that a closed, simply connected, integral, positive symplectic manifold $(X, \omega)$ has the *quasifibration property* if the evaluation map

$$E : Hol_{x_0}(\mathbb{P}^1, X)^+ \to X$$

is a quasifibration.

With these definitions, the statement of our main theorem 2 is now precise. This relates the stability condition that the holomorphic mapping space $Hol_{x_0}(\mathbb{P}^1, X)$ stabilizes to the continuous mapping space $\Omega^2 X$, to the Morse theoretic condition that the flow category of the symplectic action realizes the homotopy type of the manifold on which it is defined, $\tilde{L}X$. We will prove the theorem in the next section.

We end this section by making the following observation about how our stabilization construction is related to the group completion.

Assume that $X$ has the quasifibration property, and assume that $Hol_{x_0}(\mathbb{P}^1, X)$ has the further property that it has the structure of a $\mathcal{C}_2$ - operad space, whose $H$ - space multiplication lifts (up to homotopy) the loop sum operation of $\Omega^2 X$, and extends (up to homotopy) the monoid structure of $Hol_{x_0,x_0}(\mathbb{P}^1, X)$. The $\mathcal{C}_2$ structure assures that the monoid $\pi_0(Hol_{x_0}(\mathbb{P}^1, X))$ is commutative, and that the topological monoid $Hol_{x_0}(\mathbb{P}^1, X)$ is homotopy commutative. Assume furthermore that $\pi_0(Hol_{x_0}(\mathbb{P}^1, X))$ is finitely generated as a monoid. Let $\{\gamma_1, \cdots, \gamma_k\}$ be a set of elements in $Hol_{x_0}(\mathbb{P}^1, X)$ that generate $\pi_0(Hol_{x_0}(\mathbb{P}^1, X))$. Then in our definition of $Hol_{x_0}(\mathbb{P}^1, X)^+$ we can take our "gluing map" $\gamma = \gamma_1 + \cdots + \gamma_k$. Then by the group completion theorem [15] we have a homology equivalence

$$(1.8) \qquad Hol_{x_0}(\mathbb{P}^1, X)^+ \simeq \Omega B(Hol_{x_0}(\mathbb{P}^1, X)).$$

This structure (i.e $(X, \omega)$ a positive, integral symplectic manifold, with $Hol_{x_0}(\mathbb{P}^1, X)$ a $\mathcal{C}_2$ - operad space, with finitely generated $\pi_0$) exists, for example, when $X$ is a coadjoint orbit of a compact Lie group on its Lie algebra, or when $X$ is a loop group $X = \Omega G$, where $G$ is a simply connected compact Lie group (see [1], [23]).

## 2. Proof of theorem 2

In this section we prove the main theorem 2.

*Proof.* Throughout we will assume that $X$ has the quasifibration property. Let *Spaces* denote the category of based topological spaces and basepoint preserving continuous maps. Consider the functor

$$\mathcal{H} : \mathcal{C}_\alpha \to Spaces$$

which on objects is defined by $\mathcal{H}((x,g)) = Hol_{x_0,x}(\mathbb{P}^1, X)^+$, by which we mean the fiber at $x \in X$ of the evaluation map $E : Hol_{x_0}(\mathbb{P}^1, X)^+ \to X$. On the level of morphisms, if $\gamma \in \bar{Hol}_{x_1,x_2}(\mathbb{P}^1, X)$, then $\mathcal{H}(\gamma)$ is the gluing operation (on the right)

$$(2.1) \qquad \mathcal{H}(\gamma) : Hol_{x_0,x_1}(\mathbb{P}^1, X)^+ \xrightarrow{*\gamma} Hol_{x_0,x_2}(\mathbb{P}^1, X)^+.$$

Notice that by the quasifibration property we have the following:

**Lemma 3.** *For every morphism $\gamma \in \bar{H}ol_{x_1,x_2}(\mathbb{P}^1, X)$, the induced map*

$$\mathcal{H}(\gamma) : Hol_{x_0,x_1}(\mathbb{P}^1, X)^+ \xrightarrow{*\gamma} Hol_{x_0,x_2}(\mathbb{P}^1, X)^+.$$

*is a homotopy equivalence.*

This action of the morphisms of $\mathcal{C}_\alpha$ on the functor $\mathcal{H}$, $\mathcal{H}(x_1, g_1) \times Mor((x_1, g_1), (x_2, g_2)) \to \mathcal{H}(x_2, g_2)$, we write as

$$\mathcal{H} \times_{ob(\mathcal{C}_\alpha)} Mor(\mathcal{C}_\alpha) \to Mor(\mathcal{C}_\alpha).$$

We now consider the following "simplicial Borel construction", $\mathcal{E}_{\mathcal{C}_\alpha}(\mathcal{H})$, whose $n$ - simplices are given by

$$\mathcal{E}_{\mathcal{C}_\alpha}(\mathcal{H})_n = \mathcal{H} \times_{ob(\mathcal{C}_\alpha)} Mor(\mathcal{C}_\alpha) \times_{ob(\mathcal{C}_\alpha)} Mor(\mathcal{C}_\alpha) \times_{ob(\mathcal{C}_\alpha)} \cdots \times_{ob(\mathcal{C}_\alpha)} Mor(\mathcal{C}_\alpha),$$

$n$ - copies of $Mor(\mathcal{C}_\alpha)$. The notation $Mor(\mathcal{C}_\alpha) \times_{ob(\mathcal{C}_\alpha)} Mor(\mathcal{C}_\alpha)$ refers to taking products of composable morphisms. The face maps are defined as usual by composition of morphisms and by the action of the morphisms on the functor $\mathcal{H}$. The degeneracy maps are defined by inserting the identity morphism in the various slots. Notice that the projection maps on the level of $n$ - simplices,

$$p_n : \mathcal{E}_{\mathcal{C}_\alpha}(\mathcal{H})_n \to (B\mathcal{C}_\alpha)_n$$

$$\mathcal{H} \times_{ob(\mathcal{C}_\alpha)} Mor(\mathcal{C}_\alpha) \times_{ob(\mathcal{C}_\alpha)} \cdots \times_{ob(\mathcal{C}_\alpha)} Mor(\mathcal{C}_\alpha) \to Mor(\mathcal{C}_\alpha) \times_{ob(\mathcal{C}_\alpha)} \cdots \times_{ob(\mathcal{C}_\alpha)} Mor(\mathcal{C}_\alpha)$$

fit together to give a map of simplicial spaces

$$p : \mathcal{E}_{\mathcal{C}_\alpha}(\mathcal{H}) \to B\mathcal{C}_\alpha$$

Our main technical result needed to prove theorem 2 is the following.

**Theorem 4.** *The induced map on the level of geometric realizations,*

$$p : \|\mathcal{E}_{\mathcal{C}_\alpha}(\mathcal{H})\| \to \|B\mathcal{C}_\alpha\|$$

*is a quasifibration with fiber $\mathcal{H}(x_0) = Hol_{x_0,x_0}(\mathbb{P}^1, X)^+$. Furthermore the space $\|\mathcal{E}_{\mathcal{C}_\alpha}(\mathcal{H})\|$ is contractible.*

The following is an immediate consequence.

**Corollary 5.** *There is a natural homotopy equivalence*

$$j : Hol_{x_0,x_0}(\mathbb{P}^1, X)^+ \xrightarrow{\simeq} \Omega B\mathcal{C}_\alpha.$$

*Proof.* ( *theorem 4*). The first part of the theorem follows from lemma 3 and the following lemma, which was proven in [19].

**Lemma 6.** *If $p : E \to B$ is a map of simplicial spaces such that $E_k \to B_k$ is a quasifibration for each $k \geq 0$, and for each simplicial operation $\theta : [k] \to [l]$ and each $b \in B_l$ the map $p^{-1}(b) \to p^{-1}(\theta^*(b))$ is a homotopy equivalence, then the map of realizations $\|E\| \to \|B\|$ is a quasifibration.*

To prove the second part of the theorem we observe that the simplicial space $\mathcal{E}_{\mathcal{C}_\alpha}(\mathcal{H})$ is the nerve (classifying space) of the topological category $\mathcal{C}_\alpha(\mathcal{H})$ whose objects are elements of the space $Hol_{x_0}(\mathbb{P}^1, X)^+$, and if $\gamma_1 \in Hol_{x_0, x_1}(\mathbb{P}^1, X)^+$ and $\gamma_2 \in Hol_{x_0, x_2}(\mathbb{P}^1, X)^+$ are objects in $\mathcal{E}_{\mathcal{C}_\alpha}(\mathcal{H})$, then a morphism $\phi : \gamma_1 \to \gamma_2$ is an element of $Mor_{\mathcal{C}_\alpha}((x_1, n_1), (x_2, n_2))$ for some $n_1, n_2 \in \mathbb{Z}$, such that under the gluing operation 2.1

$$\gamma_1 * \phi = \gamma_2.$$

Furthermore, notice that the category $\mathcal{C}_\alpha(\mathcal{H})$ has an initial object: namely the constant holomorphic map $\epsilon_{x_0} \in Hol_{x_0, x_0}(\mathbb{P}^1, X)^+$. Thus the realization of its classifying space ($\cong \|\mathcal{E}_{\mathcal{C}_\alpha}(\mathcal{H})\|$) is contractible. The theorem, and hence the corollary follow. □

We remark that 4 is a kind of "group completion theorem" of the sort originally proved in [15]. Generalizations of the sort proved here (done in the category of bisimplicial sets) were done in [16], [10], and [22]. This theorem will be useful in our proof of 2, which we now complete.

Recall from section 1 the definition of the map

$$\phi : B\mathcal{C}_\alpha \to \tilde{L}X$$

via the category $\mathcal{C}(\tilde{L}X)$. This category has objects points of $\tilde{L}X$, and morphisms are paths in $\tilde{L}X$. Since an element of $\tilde{L}X$ is given by a pair $(\gamma, \theta)$, where $\gamma \in LX$, and $\theta : D^2 \to X$ is a homotopy class of an extension of $\gamma$, then a morphism $\psi$ in $\mathcal{C}(\tilde{L}X)$ from $(\gamma_0, \theta_0)$ to $(\gamma_1, \theta_1)$, can be viewed as a map of the cylinder

$$\psi : S^1 \times [0, 1] \to X$$

with $\psi_0 = \gamma_0$ and $\psi_1 = \gamma_1$. In particular, if $(\epsilon_{x_0}, \theta_0)$ and $(\epsilon_{x_1}, \theta_1)$ are the objects in $\mathcal{C}(\tilde{L}X)$ corresponding to the constant loops at $x_0$ and $x_1$ and extensions $\theta_0, \theta_1 : S^2 \to X$ with $\omega([\theta_i]) = n_i$, then a morphism in $\mathcal{C}(\tilde{L}X)$ between them are given by elements of the mapping space $Map_{x_0, x_1}^{n_0 - n_1}(S^2, X)$. Recall that the map $\phi : B\mathcal{C}_\alpha \to \tilde{L}X$ was given by the composition

$$\phi : B\mathcal{C}_\alpha \xrightarrow{j} B\mathcal{C}(\tilde{L}X) \simeq \tilde{L}X$$

where $j : \mathcal{C}_\alpha \hookrightarrow \mathcal{C}(\tilde{L}X)$ is the inclusion of categories. On morphism spaces the functor $j$ is given by including spaces of holomorphic maps into spaces of continuous maps. Now by performing the same argument used to prove 4, but replacing the category $\mathcal{C}_\alpha$ by $\mathcal{C}(\tilde{L}X)$, we see that this implies that we have the following homotopy commutative diagram:

(2.2)
$$\begin{array}{ccc} \Omega B\mathcal{C}_\alpha & \xrightarrow{\Omega\phi} & \Omega\tilde{L}X \\ \simeq \downarrow & & \downarrow \simeq \\ Hol_{x_0,x_0}(\mathbb{P}^1, X)^+ & \xrightarrow{j} & Map_{x_0,x_0}(S^2, X). \end{array}$$

Now consider the evaluation map $E : Hol_{x_0}(\mathbb{P}^1, X) \to X$. By assumption, this map is a quasifibration. In fact this means that the inclusion $j$ of holomorphic maps into continuous maps induces a map of quasifibrations

$$\begin{array}{ccc} Hol_{x_0,x_0}(\mathbb{P}^1, X)^+ & \xrightarrow{j} & Map_{x_0,x_0}(S^2, X) \\ \cap \downarrow & & \downarrow \cap \\ Hol_{x_0}(\mathbb{P}^1, X)^+ & \xrightarrow{j} & \Omega^2 X \\ E \downarrow & & \downarrow E \\ X & \xrightarrow{=} & X. \end{array}$$

Thus $j : Hol_{x_0}(\mathbb{P}^1, X)^+ \to \Omega^2 X$ is a homotopy equivalence if and only if the induced map of fibers $j : Hol_{x_0,x_0}(\mathbb{P}^1, X)^+ \to Map_{x_0,x_0}(S^2, X)$ is a homotopy equivalence. But by 2.2 this map is a homotopy equivalence if and only if $\phi : B\mathcal{C}_\alpha \to \tilde{L}X$ is a homotopy equivalence. This proves the statement of theorem 2. $\square$

## 3. Homogeneous Spaces

In this section we use theorem 2 to give an alternative proof to the following stability theorem of Gravesen [7].

**Theorem 7.** *Let $G$ be a complex linear algebraic group, and let $P < G$ be a parabolic subgroup. The homogeneous space $G/P$ has the structure of a smooth projective variety. Then $G/P$ has the stability property:*
$$Hol_{x_0}(\mathbb{P}^1, G/P)^+ \simeq \Omega^2 G/P.$$

**Remark.**

1. Gravesen stated his result in terms of the colimit under a gluing operation on $Hol_{x_0}(\mathbb{P}^1, G/P)$:
$$\varinjlim Hol_{x_0}(\mathbb{P}^1, G/P) \simeq \Omega^2 G/P.$$

    By the construction of $Hol_{x_0}(\mathbb{P}^1, G/P)^+$ in section 1, it is clear that this limit is the same as our stabilization $Hol_{x_0}(\mathbb{P}^1, G/P)^+$.

14 R.L. COHEN, J.D.S. JONES, AND G.B. SEGAL

2. In [1] Boyer, Hurtubise, Mann, and Milgram used Gravesen's theorem to prove a stronger stability theorem. Namely given $n \in \mathbb{Z}$, they showed that there is an explicit range of dimensions that increases over the limiting process, in which the inclusion

$$j_* : \pi_q(Hol^n_{x_0}(\mathbb{P}^1, G/P)) \to \pi_q(\Omega^2 G/P)$$

is an isomorphism.

By theorem 2, to prove 7 it suffices to prove the following results:

**Proposition 8.** *$G/P$ has the quasifibration property.*

**Proposition 9.** *Let $\mathcal{C}_\alpha(G/P)$ be the flow category of the symplectic action functional on $\tilde{L}(G/P)$. Then the map*

$$\phi : B\mathcal{C}_\alpha(G/P) \to \tilde{L}(G/P)$$

*is a homotopy equivalence.*

We will prove proposition 9 first.

*Proof.* For $P < G$ a parabolic subgroup, $G/P$ is a compact, smooth projective variety. Let $e : G/P \hookrightarrow \mathbb{P}^N$ be a projective embedding. The symplectic form on $G/P$ is the pull back under this embedding of the canonical symplectic form on $\mathbb{P}^N$. Now $\pi_2(\mathbb{P}^N) = \mathbb{Z}$, and the symplectic form $\omega : \pi_2(\mathbb{P}^N) \to \mathbb{Z}$ is an isomorphism. Thus $\tilde{L}\mathbb{P}^N$ is the universal cover of the loop space $L\mathbb{P}^N$. $\tilde{L}(G/P)$ is the pull back of the cover $\tilde{L}\mathbb{P}^N$ under the induced embedding of loop spaces, $e : L(G/P) \hookrightarrow L\mathbb{P}^N$.

Notice that the symplectic action functional on projective space,

$$\alpha : \tilde{L}\mathbb{P}^N \to \mathbb{R}$$

yields the symplectic action on $G/P$ via the composition

(3.1) $$\alpha : \tilde{L}(G/P) \xrightarrow{e} \tilde{L}(\mathbb{P}^N) \xrightarrow{\alpha} \mathbb{R}.$$

This induces a functor between the corresponding flow categories:

(3.2) $$\mathcal{C}_\alpha(G/P) \xrightarrow{e} \mathcal{C}_\alpha(\mathbb{P}^N).$$

In order to study the flow category $\mathcal{C}_\alpha(G/P)$ we will use this functor together with a study of $\mathcal{C}_\alpha(\mathbb{P}^N)$. This category was studied in detail in [5]. We recall those results now.

Let $W = \mathbb{C}[z, z^{-1}]$ be the vector space of Laurent polynomials, topologized as a space of maps $\mathbb{C}^\times \to \mathbb{C}$, where $\mathbb{C}^\times = \mathbb{C} - \{0\}$. The linear flow $z^k \to e^{kt}z^k$ defines a flow $\theta$ on the infinite projective space $\mathbb{P}(\mathbb{C}^{N+1} \otimes W)$. This in fact is a gradient flow. The stationary points of $\theta$ are $\mathbb{P}^N \times \mathbb{Z}$, where $\mathbb{P}^N \times \{k\}$ is the subspace $\mathbb{P}(\mathbb{C}^{N+1} \otimes z^k) \subset \mathbb{P}(\mathbb{C}^{N+1} \otimes W)$. If $W^n_m$ is the subspace of $W$ spanned by the $z^j$'s with $m \leq j \leq n$, then the space of points which lie on piecewise flow lines of $\theta$ which go from level $n$ to level $m$ is $\mathbb{P}(\mathbb{C}^{N+1} \otimes W_m^n) = \mathbb{P}^{(N+1)(n-m)-1}$, which is compact. In fact it was stressed in [5] that the flow category $\mathcal{C}_\theta$ is a compactification of $\mathcal{C}_\alpha(\mathbb{P}^N)$ in the following sense. Let $\mathcal{U}_N$ be the open dense subset of $\mathbb{P}(\mathbb{C}^{N+1} \otimes W)$ consisting of $(N+1)$ - tuples of Laurent polynomials $(p_0, \cdots, p_N) \in \mathbb{C}^{N+1} \otimes W$ with no common roots in $\mathbb{C}^\times$. (By common roots we mean roots common to $all$ the polynomials $\{p_0, \cdots, p_N\}$.) $\mathcal{U}_N$ is invariant under the flow $\theta$, and it was seen easily that the corresponding flow category is isomorphic to the category $\mathcal{C}_\alpha(\mathbb{P}^N)$. Thus we have the inclusion of flow categories, $\mathcal{C}_\alpha(\mathbb{P}^N) \subset \mathcal{C}_\theta$. Moreover, since the flow $\theta$ on $\mathbb{P}(\mathbb{C}^{N+1} \otimes W)$ is the limit of flows on the finite dimensional compact manifolds $\mathbb{P}(\mathbb{C}^{N+1} \otimes W_m^n) = \mathbb{P}^{(N+1)(n-m)-1}$, then we know that $B\mathcal{C}_\theta \cong \mathbb{P}(\mathbb{C}^{N+1} \otimes W)$. Moreover the realization of the inclusion functor $B\mathcal{C}_\alpha(\mathbb{P}^N) \subset B\mathcal{C}_\theta$ gives the inclusion of the open dense subset $\mathcal{U}_N \subset \mathbb{P}(\mathbb{C}^{N+1} \otimes W)$.

In particular if $\mathcal{C}_\theta^{n,m}$ is the full subcategory of $\mathcal{C}_\theta$ whose objects are $\mathbb{P}^N \times \{m, m+1, \cdots, n\}$, then by the results of [4] the map

$$\phi : B\mathcal{C}_\theta^{n,m} \to \mathbb{P}(\mathbb{C}^{N+1} \times W_m^n)$$

is a homotopy equivalence, where the inverse image of each point is a simplex. Taking the limit over $n$ and $m$, we have that

(3.3) $$\phi : B\mathcal{C}_\theta \to \mathbb{P}(\mathbb{C}^{N+1} \times W)$$

is a proper map and a homotopy equivalence, where the inverse image of each point is a simplex. The pull back of $\phi$ to $\mathcal{U}_N \subset \mathbb{P}(\mathbb{C}^{N+1} \times W)$ is the map $B\mathcal{C}_\alpha \to \mathcal{U}_N$, which is therefore a homotopy equivalence.

Now it was seen in [5] that the map $\phi : B\mathcal{C}_\alpha(\mathbb{P}^N) \to \tilde{L}(\mathbb{P}^N)$ can be realized in terms of the space $\mathcal{U}_N$, by observing that an $(N+1)$ - tuple of Laurent polynomials

$$(p_0, \cdots, p_N) \in \mathcal{U}_N \subset \mathbb{P}(\mathbb{C}^{N+1} \otimes W)$$

determines a map

$$p : \mathbb{C}^\times \to \mathbb{P}^N.$$

Since this map is algebraic it extends to a holomorphic map $p : \mathbb{C} \cup \infty \to \mathbb{P}^N$. By restricting this holomorphic map to the equator (i.e the unit circle $S^1 \subset \mathbb{C}$) we get an element of $L\mathbb{P}^N$. Using the given extension of this map to the unit disk, this actually defines an element of $\tilde{L}\mathbb{P}^N$. This defines an embedding $\mathcal{U}_N \hookrightarrow \tilde{L}\mathbb{P}^N$, which makes the following diagram homotopy commute:

(3.4)
$$\begin{array}{ccc} B\mathcal{C}_\alpha(\mathbb{P}^N) & \xrightarrow{\phi} & \tilde{L}\mathbb{P}^N \\ \simeq \downarrow & & \downarrow = \\ \mathcal{U}_N & \hookrightarrow & \tilde{L}\mathbb{P}^N. \end{array}$$

This describes the homotopy type of the classifying space $B\mathcal{C}_\alpha(\mathbb{P}^N)$ as realizing the subspace of $\tilde{L}(\mathbb{P}^N)$ consisting of those loops whose projective coordinates have finite Fourier expansions (i.e are Laurent polynomials). We see that by restricting to $\mathcal{C}_\alpha(G/P)$, we get a similar description for $B\mathcal{C}_\alpha(G/P)$ as follows.

The projective embedding $e : G/P \hookrightarrow \mathbb{P}^N$ is defined by a homogeneous ideal $I(G/P) \subset \mathbb{C}[x_0, \cdots, x_n]$. Then consider the subspace $\mathcal{U}(G/P) \subset \mathcal{U}_N \subset \mathbb{P}(\mathbb{C}^{N+1} \otimes W)$ defined by

(3.5) $\quad \mathcal{U}(G/P) = \{(p_0, \cdots, p_N) \in \mathcal{U}_N : f(p_0, \cdots, p_N) = 0 \,\text{for all}\, f \in I(G/P)\}.$

Notice that $\mathcal{U}(G/P) \subset \mathcal{U}_N$ is a $\theta$ - flow invariant subspace. Its stationary points correspond to $G/P \times \mathbb{Z} \subset \mathbb{P}^n \times \mathbb{Z}$, and the points lying on flows going from level $n$ to level $m$ are given by $(N+1)$ - tuples, $(p_0, \cdots, p_N) \subset \mathbb{P}(C^{N+1} \otimes W_m^n)$ with $f(p_0, \cdots, p_n) = 0$ for all $f \in I(G/P)$. This space precisely parameterizes the holomorphic maps $\gamma : \mathbb{P}^1 \to G/P$ so that the composition $\mathbb{P}^1 \to G/P \hookrightarrow \mathbb{P}^N$ has degree $n - m$. Thus the flow category of $\mathcal{U}(G/P)$ is isomorphic to $\mathcal{C}_\alpha(G/P)$. In particular the pull back of the map $\phi : B\mathcal{C}_\theta \to \mathbb{P}(\mathbb{C}^{N+1} \times W)$ to the subspace $\mathcal{U}(G/P) \subset \mathcal{U}_N \subset \mathbb{P}(\mathbb{C}^{N+1} \times W)$ is precisely $B\mathcal{C}_\alpha(G/P)$. Since $\phi$ is a proper map and a homotopy equivalence with the inverse image of each point a simplex (3.3), then the pull back $B\mathcal{C}_\alpha(G/P) \to \mathcal{U}(G/P)$ is a homotopy equivalence. Note also that the image of $\mathcal{U}(G/P) \subset \mathcal{U}_N \hookrightarrow \tilde{L}\mathbb{P}^N$ lies in $\tilde{L}(G/P)$ and consists precisely of those loops whose homogeneous coordinates have finite Laurent expansion. So we have the following commutative diagram:

(3.6)
$$\begin{array}{ccc} B\mathcal{C}_\alpha(G/P) & \xrightarrow{e} & B\mathcal{C}_\alpha(\mathbb{P}^N) \\ \simeq \downarrow & & \downarrow \simeq \\ \mathcal{U}(G/P) & \hookrightarrow & \mathcal{U}_N \\ \cap \downarrow & & \downarrow \cap \\ \tilde{L}(G/P) & \xrightarrow{e} & \tilde{L}(\mathbb{P}^N) \end{array}$$

Moreover the vertical compositions in this diagram are homotopic to the maps $\phi : B\mathcal{C}_\alpha(G/P) \to \tilde{L}(G/P)$ and $\phi : B\mathcal{C}_\alpha(\mathbb{P}^N) \to \tilde{L}(\mathbb{P}^N)$. Thus to prove that $\phi : B\mathcal{C}_\alpha(G/P) \to \tilde{L}(G/P)$ is a homotopy equivalence it suffices to prove the following.

**Lemma 10.** *The inclusion $\mathcal{U}(G/P) \hookrightarrow \tilde{L}(G/P)$ is a homotopy equivalence.*

**Remark.** As observed above, this inclusion can be viewed as the inclusion of the space of *polynomial loops* (i.e those loops which, in homogeneous coordinates have finite Fourier expansion) into the space of all loops.

*Proof.* Since $G$ is a linear algebraic group, it is a subgroup of $GL(N, \mathbb{C})$ for some $N$. This $N$ can be taken to be the dimension of the projective embedding $e : G/P \hookrightarrow \mathbb{P}^N$. So in particular a loop in $G$ can be viewed as a loop in the affine space of matrices. Let $LG$ be the infinite group of smooth loops in $G$, and $L_{pol}G$ the subgroup of *polynomial* loops; those loops in $G$ which, together with their inverses have finite Fourier expansion. That is, they are given by finite Laurent polynomials. Notice that $L(GL(N, \mathbb{C}))$ acts transitively on $\tilde{L}(\mathbb{P}^N)$, and the action restricts to a transitive action of $LG$ on $\tilde{L}(G/P)$. Furthermore the isotropy subgroup of a constant loop is clearly a union of path components of the loop group of $P$, which we denote by $L(P)_0$. This subgroup of $L(P)$ has the property that the quotient group is infinite cyclic and that the projection map $p : L(G)/L(P)_0 \to L(G)/L(P) = L(G/P)$ is the infinite cyclic cover $p : \tilde{L}(G/P) \to L(G/P)$.

Now the polynomial loop group $L_{pol}(GL(N, \mathbb{C}))$ acts transitively on the space of polynomial loops $\mathcal{U}_N \subset L(\mathbb{P}^N)$. It restricts to give a transitive action of $L_{pol}(G)$ on $\mathcal{U}(G/P)$. In this case the isotropy subgroup of a constant loop is th union of path components $L_{pol}(P)_0$ of $L_{pol}(P)$ corresponding to the subgroup $L(P)_0$ of $L(P)$. We have therefore shown that there is a homeomorphism $L_{pol}G/L_{pol}(P)_0 \cong \mathcal{U}(G/P)$ making the following diagram commute:

(3.7)
$$\begin{array}{ccc} L_{pol}G/L_{pol}(P)_0 & \xrightarrow{\cong} & \mathcal{U}(G/P) \\ \cap \downarrow & & \downarrow \cap \\ L(G)/L(P)_0 & \xrightarrow{\cong} & \tilde{L}(G/P). \end{array}$$

Now it was proved in chapter 8 of [17] that the inclusion of the polynomial loop group into the smooth loop group $L_{pol}(G) \hookrightarrow L(G)$ is a homotopy equivalence. (This proof also uses Morse theory!) Since this also holds for $L_{pol}(P) \hookrightarrow L(P)$ it also holds for the corresponding union of connected components $L_{pol}(P)_0 \xrightarrow{\simeq} L(P)_0$. This implies that the inclusion $L_{pol}G/L_{pol}(P)_0 \hookrightarrow L(G)/L(P)_0$ is a homotopy equivalence. By by this diagram this implies that the inclusion $\mathcal{U}(G/P) \hookrightarrow \tilde{L}(G/P)$ is a homotopy equivalence. This proves the lemma.

□

Now as remarked earlier, the proof of the lemma completes the proof of proposition 9.  □

We now prove proposition 8. As seen earlier, this is the last step in the proof of theorem 7.

*Proof.* We need to show that for $X = G/P$, the evaluation map
$$E : Hol_{x_0}(\mathbb{P}^1, X)^+ \to X$$
is a quasifibration. Let $x_1 \in X$. It is well known that the homogeneous space $X = G/P$ is rationally connected (see for example [12]) so $Hol_{x_0, x_1}(\mathbb{P}^1, X)$ is nonempty. Let $\theta \in Hol_{x_0, x_1}(\mathbb{P}^1, X)$ be any element. Gluing with $\theta$ induces a map $\theta_* : Hol_{x_0, x_0}(\mathbb{P}^1, X)^+ \to Hol_{x_0, x_1}(\mathbb{P}^1, X)^+$. It suffices to show that
$$\theta_* : Hol_{x_0, x_0}(\mathbb{P}^1, X)^+ \to Hol_{x_0, x_1}(\mathbb{P}^1, X)^+$$
is a homotopy equivalence. Now since $X = G/P$ is rationally connected, $Hol_{x_1, x_0}(\mathbb{P}^1, X)$ is nonempty. Let $\alpha \in Hol_{x_1, x_0}(\mathbb{P}^1, X)$. We will prove that gluing with the products,
$$(\alpha * \theta)_* : Hol_{x_0, x_0}(\mathbb{P}^1, X)^+ \to Hol_{x_0, x_0}(\mathbb{P}^1, X)^+$$

and

$$(\theta * \alpha)_* : Hol_{x_0,x_1}(\mathbb{P}^1, X)^+ \to Hol_{x_0,x_1}(\mathbb{P}^1, X)^+$$

are homotopy equivalences. Now for $X = G/P$, the holomorphic mapping space $Hol_{x_0}(\mathbb{P}^1, X)$ is a $\mathcal{C}_2$ - operad space. The structure is studied, for example, in [1]. This in particular implies that the monoid structure in $Hol_{x_0,x_0}(\mathbb{P}^1, X)$ is homotopy commutative. Furthermore $\pi_0(Hol_{x_0}(\mathbb{P}^1, X))$ is finitely generated [1], [7]. So by 1.8 this implies that $Hol_{x_0}(\mathbb{P}^1, X)^+$ is group complete, and so in particular $\pi_0(Hol_{x_0}(\mathbb{P}^1, X)^+)$ is a group. But since $X$ is rationally connected, $\pi_0(Hol_{x_0,x_0}(\mathbb{P}^1, X)^+) \cong \pi_0(Hol_{x_0}(\mathbb{P}^1, X)^+)$. Thus $\pi_0(Hol_{x_0,x_0}(\mathbb{P}^1, X)^+)$ is a group. This means that the class $[\alpha * \theta] \in \pi_0(Hol_{x_0,x_0}(\mathbb{P}^1, X)^+)$ has an inverse, and so the element $\alpha * \theta \in Hol_{x_0,x_0}(\mathbb{P}^1, X)^+$ has a homotopy inverse. This means that $(\alpha*\theta)_* : Hol_{x_0,x_0}(\mathbb{P}^1, X)^+ \to Hol_{x_0,x_0}(\mathbb{P}^1, X)^+$ is a homotopy equivalence.

Now $X = G/P$ is homogeneous, so $x_1 = gx_0$ for some $g \in G$. Multiplication by $g$ is a holomorphic map from $G/P$ to itself, and so the class $\theta * \alpha \in Hol_{x_1,x_1}(\mathbb{P}^1, X)$ determines an element $g \cdot (\theta * \alpha) \in Hol_{x_0,x_0}(\mathbb{P}^1, X)$. As argued above, this class has a homotopy inverse which we call $\beta \in Hol_{x_0,x_0}(\mathbb{P}^1, X)^+$. Notice then that $g^{-1} \cdot \beta \in Hol_{x_1,x_1}(\mathbb{P}^1, X)^+$ defines a gluing map

$$(g^{-1} \cdot \beta)_* : Hol_{x_0,x_1}(\mathbb{P}^1, X)^+ \to Hol_{x_0,x_1}(\mathbb{P}^1, X)^+$$

which is a homotopy inverse of $(\theta*\alpha)_*$. Thus $(\theta*\alpha)_*$ is a homotopy equivalence as well. This proves the proposition and thereby completes the proof of theorem 7. □

DEPT. OF MATHEMATICS, STANFORD UNIVERSITY, STANFORD, CALIFORNIA 94305

*E-mail address*, Cohen: `ralph@math.stanford.edu`

DEPARTMENT OF MATHEMATICS, UNIVERSITY OF WARWICK, COVENTRY, ENGLAND

*E-mail address*, Jones: `jdsj@maths.warwick.ac.uk`

DEPT. OF PURE MATHEMATICS AND MATHEMATICAL STATISTICS, CAMBRIDGE UNIVERSITY, CAMBRIDGE, ENGLAND

*E-mail address*, Segal: `G.B.Segal@dpmms.cam.ac.uk`